\documentclass{article}

\usepackage[all]{xy}
\usepackage{latexsym,graphics,epsfig}
\usepackage{rotate}
\usepackage{amsfonts,amssymb}

\title{Higher Gauge Theory}

\author{John Baez
\\ Department of Mathematics
\\ University of California
\\ Riverside, CA 92521, USA 
\\
\\ Urs Schreiber 
\\ Fachbereich Mathematik
\\ Universit{\"a}t Hamburg
\\ Hamburg, 20146, Germany
\\
\\ {\small email: baez@math.ucr.edu,
    urs.schreiber@math.uni-hamburg.de}
\\
\\ {\small Hamburg Preprint ZMP-HH/05-25}
}

\date{{\small November 28, 2005}}

\newcommand{\R}{{\mathbb R}}

\newcommand{\Z}{{\mathbb Z}}

\newcommand{\maps}{\colon}

\def\stackto #1 { \, {\stackrel{#1}{\longrightarrow}}\, }
\def\stackTo #1 { {\stackrel{#1}{\Longrightarrow}} }
\newcommand{\To}{\Rightarrow}

\newcommand{\U}{{\rm U}}

\newcommand{\Aut}{{\rm Aut}}
\newcommand{\AUT}{\mathcal{AUT}}

\newcommand{\g}{\mathfrak{g}}
\newcommand{\h}{\mathfrak{h}}

\def\semidir {\ltimes}

\newcommand{\Cat}{\mathrm{Cat}}
\newcommand{\Grp}{\mathrm{Grp}}
\newcommand{\Lie}{\mathrm{Lie}}
\newcommand{\Alg}{\mathrm{Alg}}
\newcommand{\Diff}{\mathrm{Diff}}

\newcommand{\Ob}{\mathrm{Ob}}
\newcommand{\Mor}{\mathrm{Mor}}

\renewcommand{\P}{\mathcal{P}}

\def\hol {{\rm hol}}

   \def\twogroup {\mathcal{G}}

   \def\G {G}

   \def\manifold {{M}}

   \def\twoF {{F}}

   \def\twoU {{U}}

\newtheorem{theorem}{Theorem}    

\newtheorem{definition}[theorem]{Definition}
\newtheorem{example}[theorem]{Example}

\newcommand{\et}{\hspace{-0.08in}{\bf .}\hspace{0.1in}}
\def\endofproof {\hfill$\Box$\newline}
\newcommand\Proof{\noindent\emph{Proof.}\ }

\def\refdef #1{{Def.\ \ref{#1}}}

\def\of #1{\!\left({#1}\right)}
\def\id {i}
\def\extd {\mathbf{{d}}}

\def\iso{\cong}

\begin{document}
\maketitle

\begin{abstract}
\noindent
  Just as gauge theory describes the parallel transport of point
  particles using connections on bundles, higher gauge theory 
  describes the parallel transport of 1-dimensional objects (e.g.\
  strings) using 2-connections on 2-bundles.  A 2-bundle is a 
  categorified version of a bundle: that is, one where the fiber
  is not a manifold but a category with a suitable smooth structure.
  Where gauge theory uses Lie groups and Lie algebras, higher gauge
  theory uses their categorified analogues: Lie 2-groups and Lie 2-algebras.
  We describe a theory of 2-connections on principal 2-bundles and 
  explain how this is related to Breen and Messing's theory of connections
  on nonabelian gerbes.  The distinctive feature of our theory is that
  a 2-connection allows parallel transport along paths and surfaces 
  in a parametrization-independent way.  In terms of Breen and Messing's 
  framework, this requires that the `fake curvature' must vanish.
  In this paper we summarize the main results of our theory without
  proofs.
\end{abstract}


\section{Introduction}
\label{section: Introduction}

Ordinary gauge theory describes how 0-dimensional particles
transform as we move them along 1-dimensional
paths.  It is natural to assign a group element to 
each path:
\[
\xymatrix{
   \bullet \ar@/^1pc/[rr]^{g}
&& \bullet
}
\]
The reason is that composition of paths then corresponds to
multiplication in the group:
\[
\xymatrix{
   \bullet \ar@/^1pc/[rr]^{g}
&& \bullet \ar@/^1pc/[rr]^{g'}
&& \bullet
}
\]
while reversing the direction of a path corresponds to taking
inverses:
\[
\xymatrix{
   \bullet
&& \bullet \ar@/_1pc/[ll]_{g^{-1}}
}
\]
and the associative law makes the holonomy along a triple
composite unambiguous:
\[
\xymatrix{
   \bullet \ar@/^1pc/[rr]^{g}
&& \bullet \ar@/^1pc/[rr]^{g'}
&& \bullet \ar@/^1pc/[rr]^{g''}
&& \bullet
}
\]
In short, the topology dictates the algebra!

Now suppose we wish to do something similar
for 1-dimensional `strings' that trace out 2-dimensional surfaces
as they move.  Naively we might wish our holonomy to assign a group 
element to each surface like this:
\[
\xymatrix{
   \bullet \ar@/^1pc/[rr]_{}="0"
           \ar@/_1pc/[rr]_{}="1"
           \ar@{=>}"0";"1"^{g}
&& \bullet
}
\]
There are two obvious ways to compose surfaces of this sort, vertically:
\[
\xymatrix{
   \bullet \ar@/^2pc/[rr]_{}="0"
           \ar[rr]_{}="1"
           \ar@{=>}"0";"1"^{g}
           \ar@/_2pc/[rr]_{}="2"
           \ar@{=>}"1";"2"^{g'}
&& \bullet
}
\]
and horizontally:
\[
\xymatrix{
   \bullet \ar@/^1pc/[rr]_{}="0"
           \ar@/_1pc/[rr]_{}="1"
           \ar@{=>}"0";"1"^{g}
&& \bullet \ar@/^1pc/[rr]_{}="2"
           \ar@/_1pc/[rr]_{}="3"
           \ar@{=>}"2";"3"^{g'}
&& \bullet
}
\]
Suppose that both of these correspond to multiplication
in the group $G$.  Then to obtain a well-defined holonomy for this
surface regardless of whether we do vertical or horizontal composition
first:
\[
\xymatrix{
   \bullet \ar@/^2pc/[rr]_{}="0"
           \ar[rr]_{}="1"
           \ar@{=>}"0";"1"^{g_1}
           \ar@/_2pc/[rr]_{}="2"
           \ar@{=>}"1";"2"^{g_1^\prime}
&& \bullet \ar@/^2pc/[rr]_{}="3"
           \ar[rr]_{}="4"
           \ar@{=>}"3";"4"^{g_2}
           \ar@/_2pc/[rr]_{}="5"
           \ar@{=>}"4";"5"^{g'_2}
&& \bullet
}
\]
we must have
\[      (g_1 g_2)(g_1' g_2') = (g_1 g_1')(g_2 g_2') . \]
This forces $G$ to be abelian!

In fact, this argument goes back to a classic paper by
Eckmann and Hilton \cite{EckmannHilton}.  They showed
that even if we allow $G$ to be equipped with two products,
say $gg'$ for vertical composition and $g \circ g'$ for
horizontal, so long as both products share the
same unit and satisfy this `{\bf interchange law}':
\[  
  (g_1 g_1') \circ (g_2 g_2') =
      (g_1 \circ g_2) (g_1' \circ g_2')  
\]
then in fact they must agree --- so by the previous argument,
both are abelian.  The proof is very easy:
\[  
     g \circ g' 
   = (g 1) \circ (1 g') 
   = (g \circ 1) (1 \circ g') 
   = g g' 
\]

Pursuing this approach, we would ultimately reach the theory of
connections on abelian gerbes [2-8].
If $G = \U(1)$, such a connection can be locally identified with a 
2-form --- but globally it is a subtler object, just as a connection
on a $\U(1)$ bundle can be locally identified with a 1-form, but not
globally.  In fact, connections on abelian gerbes play an important role 
in string theory [9-11].
Just as ordinary electromagnetism is described by a connection on a $\U(1)$ 
bundle, usually called the `vector potential' and denoted $A$, the 
stringy analogue of electromagnetism is described by a connection on 
a $\U(1)$ gerbe, called the $B$ field.

To go beyond this and develop a theory of \emph{nonabelian} higher gauge 
fields, we must let the topology dictate the algebra.  Readers familiar
with higher categories will already have noticed that 1-dimensional 
pictures above resemble diagrams in category theory, while the
2-dimensional pictures resemble diagrams in 2-category theory.  This 
suggests that the holonomies in higher gauge 
theory should take values in some `categorified' analogue of
a Lie group --- that is, some gadget resembling a Lie group, but 
which is a category rather than a set.  We call this `Lie 2-group'.   

In fact, Lie 2-groups and their Lie 2-algebras have already been studied 
\cite{BaezLauda:2003,BaezCrans:2003} and interesting examples have been 
constructed using the mathematics of string theory: central extensions 
of loop groups \cite{BCSS}.  
But even without knowing this, we could be led to the definition of a 
Lie 2-group by considering a kind of connection that gives holonomies 
\emph{both for paths and for surfaces}.

So, let us assume that for each path we have a holonomy taking values 
in some Lie group $G$, where composition of paths corresponds to 
multiplication in $G$.  Assume also that for each 1-parameter family of 
paths with fixed endpoints we have a holonomy taking values in some other
Lie group $H$, where vertical composition corresponds to multiplication in
$H$:
\[
\xymatrix{
   \bullet \ar@/^2pc/[rr]_{}="0"
           \ar[rr]_{}="1"
           \ar@{=>}"0";"1"^{h}
           \ar@/_2pc/[rr]_{}="2"
           \ar@{=>}"1";"2"^{h'}
&& \bullet
}
\]

Next, assume that we can parallel transport an element $g \in G$
along a 1-parameter family of paths to get a new element $g' \in G$:
\[
\xymatrix{
   \bullet \ar@/^1pc/[rr]^{g}_{}="0"
           \ar@/_1pc/[rr]_{g'}_{}="1"
           \ar@{=>}"0";"1"^{h}
&& \bullet
}
\]
This picture suggests that we should think of $h$ as a kind of
`arrow' or `morphism' going from $g$ to $g'$.  We can use categories
to formalize this.  In category theory, when a morphism
goes from an object $x$ to an object $y$, we think of the morphism as
determining both its source $x$ and its target $y$.  The group
element $h$ does not determine $g$ or $g'$.  However, the pair $(g,h)$
does.  

For this reason, it is useful to create a category $\twogroup$ where the
set of objects, say $\Ob(\twogroup)$, is just $G$, while the set of
morphisms, say $\Mor(\twogroup)$, consists of ordered pairs $f = (g,h) 
\in G \times H$.  Switching our notation to reflect this, we rewrite the 
above picture as
\[
\xymatrix{
   \bullet \ar@/^1pc/[rr]^{g}_{}="0"
           \ar@/_1pc/[rr]_{g'}_{}="1"
           \ar@{=>}"0";"1"^{f}
&& \bullet
}
\]
and write $f \maps g \to g'$ for short.  

In this new notation, we can vertically compose $f \maps g \to g'$ 
and $f' \maps g' \to g''$ to get $f f' \maps g \to g''$, as follows:
\[
\xymatrix{
   \bullet \ar@/^2pc/[rr]^{g}_{}="0"
           \ar[rr]^<<<<<<{g'}_{}="1"
           \ar@{=>}"0";"1"^{f}
           \ar@/_2pc/[rr]_{g''}_{}="2"
           \ar@{=>}"1";"2"^{f'}
&& \bullet
}
\]
This is just composition of morphisms in the category $\twogroup$.  However, 
we can also horizontally compose $f_1 \maps g_1 \to g_1'$ and $f_2 \maps
g_2 \to g_2'$ to get $f_1 \circ f_2 \maps g_1g_2 \to g_1'g_2'$, as follows:
\[
\xymatrix{
   \bullet \ar@/^1pc/[rr]^{g_1}_{}="0"
           \ar@/_1pc/[rr]_{g_1'}_{}="1"
           \ar@{=>}"0";"1"_{f_1}
&& \bullet \ar@/^1pc/[rr]^{g_2}_{}="2"
           \ar@/_1pc/[rr]_{g_2'}_{}="3"
           \ar@{=>}"2";"3"_{f_2}
&& \bullet
}
\]
We assume this operation makes $\Mor(\twogroup)$ into a 
group with the pair $(1,1) \in G \times H$ as its multiplicative unit.

The good news is that now we can assume an interchange law saying this
holonomy is well-defined:
\[
\xymatrix{
   \bullet \ar@/^2pc/[rr]^{g_1}_{}="0"
           \ar[rr]^<<<<<<{g_2}_{}="1"
           \ar@{=>}"0";"1"^{f_1}
           \ar@/_2pc/[rr]_{g_3}_{}="2"
           \ar@{=>}"1";"2"^{f_1'}
&& \bullet \ar@/^2pc/[rr]^{g_1'}_{}="3"
           \ar[rr]^<<<<<<{g_2'}_{}="4"
           \ar@{=>}"3";"4"^{f_2}
           \ar@/_2pc/[rr]_{g_3'}_{}="5"
           \ar@{=>}"4";"5"^{f_2'}
&& \bullet
}
\]
namely:
\[
  (f_1 f_1') \circ (f_2 f_2') =  (f_1 \circ f_2) (f_1' \circ f_2')
\]
without forcing either $G$ or $H$ to be abelian!   Instead, the group
$\Mor(\twogroup)$ is forced to be a semidirect product of $G$ and $H$.

The structure we are rather roughly describing here turns out to be
none other than a `Lie 2-group'.  This is an `internal category' in the
category of Lie groups.  In other words, it is a category where the set of 
objects is a Lie group, the set of morphisms is a Lie group, and all the 
usual category operations are Lie group homomorphisms. 

This audacious process --- taking a familiar 
mathematical concept defined using sets and function 
and transplanting it to live within some other category --- is far from new.  
For example, a Lie group is a group in $\Diff$, the category of smooth 
manifolds and smooth maps.  The general idea of an `internal category'
living within some other category was described by Charles Ehresmann 
\cite{Ehresmann:1963} in the early 1960's.  In the next section we begin 
by reviewing this idea.  In the rest of the paper, we use it to categorify 
the theory of Lie groups, Lie algebras, bundles and connections.  
Before we proceed, let us sketch our overall plan.

The starting point is the ordinary concept of a 
principal fiber bundle.  Such a bundle can be specified using the 
following `gluing data': 
\begin{itemize}
  \item
    a base manifold $M$,
  \item
    a cover of $M$ by open sets $\{U_i\}_{i\in I}$,
  \item
    a Lie group $G$ (the `gauge group' or `structure group'), 
  \item
    on each double overlap $U_{ij} = U_i \cap U_j$ a 
    $G$-valued function $g_{ij}$, 
  \item
    such that on triple overlaps
    the following transition law holds:
    \begin{eqnarray}
       g_{ij} g_{jk} = g_{ik}.
       \nonumber
    \end{eqnarray}
\end{itemize}
Such a bundle is augmented with a connection by specifying:
\begin{itemize}
  \item
    in each open set $U_i$ a smooth functor $\mathrm{hol}_i \maps
    \P_1(U_i) \to G$ from the path groupoid of $U_i$ to the gauge group, 
  \item
   such that for all paths $\gamma$ in double overlaps $U_{ij}$
   the following transition law holds:
   \begin{eqnarray}
     \mathrm{hol}_i\of{\gamma} &=& 
      g_{ij} \, \mathrm{hol}_j\of{\gamma} \,  g_{ij}^{-1}.
     \nonumber
   \end{eqnarray}
\end{itemize}
Here the `path groupoid' $\P_1(\manifold)$ of a manifold $\manifold$ 
has points of $\manifold$ as objects and certain equivalence classes 
of smooth paths in $\manifold$ as morphisms.  There are various ways 
to work out the technical details and make $\P_1(\manifold)$ into a 
`smooth groupoid'; here we follow the approach of Barrett \cite{Barrett:1991}, 
who uses `thin homotopy classes' of paths.  Technical details aside, 
the basic idea is that a connection on a trivial $G$-bundle gives a 
well-behaved map assigning to each path $\gamma$ in the base space the 
holonomy $\hol(\gamma) \in G$ of the connection along that path.   
Saying this map is a `smooth functor' means that these holonomies 
compose when we compose paths, and that the holonomy $\hol(\gamma)$ 
depends smoothly on the path $\gamma$.

A basic goal of higher gauge theory is to categorify all of this and to work 
out the consequences.  As mentioned, the key tool is internalization.  
This leads us immediately to the concept of a Lie 2-group, and also 
to that of a `smooth 2-space': a category in $\Diff$, or more generally 
in some category of smooth spaces that allows for infinite-dimensional 
examples.  

Using these concepts, Bartels \cite{Bartels:2004} has defined 
a `principal 2-bundle $E$ over $M$ with structure 
2-group $\twogroup$'.  To arrive at this definition, the 
key steps are to replace the total space $E$ and base space $M$ of
a principal bundle by smooth 2-spaces, and to replace the structure 
group by a Lie 2-group.  In this paper we only consider the case where 
$M$ is an ordinary space, which can be regarded as a 2-space with only 
identity morphisms.  We show that for a suitable choice of structure 
2-group, principal 2-bundles give abelian gerbes over $M$.  For another
choice, they give nonabelian gerbes.  This sets the stage for a result 
relating the 2-bundle approach to higher gauge theory to Breen and 
Messing's approach based on nonabelian gerbes \cite{BreenMessing:2001}.

Just as a connection on a trivial principal bundle over $M$ gives a 
functor from the path groupoid of $M$ to the structure group, one 
might hope that a `2-connection' on a trivial principal 2-bundle would 
define a 2-functor from some sort of `path 2-groupoid' to the structure 
2-group.  This has already been confirmed in the context of higher lattice 
gauge theory [19-21].
Thus, the main issues not yet addressed are those involving differentiability. 

To address these issues, we define for any smooth space $M$ a smooth 
2-groupoid $\P_2(M)$ such that:
\begin{itemize}
\item
the objects of $\P_2(M)$ are points of $M$: 
$ \quad   \bullet \, {\textstyle{\small{\textit{x}}}}   $
\item 
the morphisms of $\P_2(M)$ are thin homotopy classes of
smooth paths $\gamma \maps [0,1] \to M$ such that $\gamma(s)$ 
is constant in a neighborhood of $s = 0$ and $s = 1$:
$
\xymatrix{
 {\textstyle{\small{\textit{x}}}}  \,
   \bullet \ar@/^1pc/[rr]^{\gamma}
&& \bullet \, {\textstyle{\small{\textit{y}}}}   
}
$
\item 
the 2-morphisms of $\P_2(M)$ are `bigons': that is,
thin homotopy classes of smooth maps 
$\Sigma \maps [0,1]^2 \to M$ such that $\Sigma(s,t)$ is 
constant near $s = 0$ and $s = 1$, and 
independent of $t$ near $t = 0$ and $t = 1$:
$
\xymatrix{
 {\textstyle{\small{\textit{x}}}} \, 
   \bullet \ar@/^1pc/[rr]^{\gamma_1}_{}="0"
           \ar@/_1pc/[rr]_{\gamma_2}^{}="1"
&& \bullet \, {\textstyle{\small{\textit{y}}}}   
\ar@{=>}"0";"1"^{\Sigma}
}
$
\end{itemize}
The `thin homotopy' equivalence relation, borrowed from the work
of Mackaay and Picken \cite{MackaayPicken:2000,Mackaay:2001}, 
guarantees that two maps differing only by a reparametrization 
define the same bigon.  This is important because we seek a 
{\it reparametrization-invariant notion of surface holonomy}.

We define a `2-connection' on a trivial principal 2-bundle 
over $M$ to be a smooth 2-functor $\hol \maps \P_2(M) 
\to \twogroup$, where $\twogroup$ is the structure 2-group.   
This means that the 2-connection assigns holonomies both
to paths and bigons, independent of their parametrization, compatible 
with the standard operations of composing paths and bigons, and 
depending smoothly on the path or surface in question.  

We also define 2-connections for nontrivial principal 2-bundles,
and state a theorem obtaining these from Lie-algebra-valued
differential forms.   We then show that for a certain class of
structure 2-groups, such differential forms reduce to Breen and 
Messing's `connections on nonabelian gerbes' \cite{BreenMessing:2001}.  
The surprise is that we only obtain connections satisfying a 
certain constraint: the `fake curvature' must vanish!

To understand this, one must recall \cite{BaezLauda:2003}
that a Lie 2-group $\twogroup$ amounts to the same thing as 
a `crossed module' of Lie groups $(G,H,t,\alpha)$, where:
\begin{itemize}
\item $G$ is the group of objects of $\G$, $\Ob(\twogroup)$:
\item $H$ is the subgroup of $\Mor(\twogroup)$ consisting
of morphisms with source equal to $1 \in G$:
\item $t \maps H \to G$ is the homomorphism sending each morphism
in $H$ to its target,
\item $\alpha$ is the action of $G$ as automorphisms of $H$ defined
using conjugation in $\Mor(\twogroup)$ as follows: 
$\alpha(g) h = 1_g h {1_g}^{-1}$. 
\end{itemize}
Differentiating all this data one obtains a `differential crossed 
module' $(\g,\h,dt,d\alpha)$, which is just another way of talking
about a Lie 2-algebra \cite{BaezSchreiber:2004}.

In these terms, a 2-connection on a trivial principal 2-bundle over
$M$ with structure 2-group $\twogroup$ consists of a $\g$-valued
1-form $A$ together with an $\h$-valued 2-form $B$ on $M$. Translated
into this framework, Breen and Messing's `fake curvature' is the
$\g$-valued 2-form
\[
dt\of{B} + F_A, 
\]
where $F_A = \extd A + A \wedge A$ is the usual curvature of $A$.
We show that {\it if and only if the fake curvature vanishes}, 
one obtains a well-defined 2-connection $\hol \maps \P_2(M) \to \twogroup$.  

The importance of vanishing fake curvature in the framework of lattice 
gauge theory was already emphasized by Girelli and Pfeiffer
\cite{GirelliPfeiffer:2004}.  The special case where also $F_A = 0$ 
was studied by Alvarez, Ferreira, Sanchez and Guillen
\cite{AlvarezFerreiraSanchezGuillen:1998}.   The case where 
$G = H$ has been studied already by the second author of this
paper \cite{Schreiber:2004e}.  Our result subsumes these 
cases in a common framework.  

This paper is an introduction to work in progress \cite{BaezSchreiber:2004},
which began in rudimentary form as an article by the first author
\cite{Baez:2002}, and overlaps to some extent with theses by 
Bartels \cite{Bartels:2004} and the second author
\cite{Schreiber:2005}.  Bartels' thesis develops the general
theory of 2-bundles.  The second author's thesis investigates the
relationship between nonabelian higher gauge theory and the physics 
of strings.  Aschieri, Cantini and Jur{\v c}o 
\cite{AschieriCantiniJurco:2003,AschieriJurco:2004} have
also studied this subject, using connections on nonabelian gerbes.
Other physicists, including Chepelev \cite{Chepelev:2002} and Hofman 
\cite{Hofman:2002}, have also studied nonabelian higher gauge fields.

\section{Internalization}
\label{section: Internalization}

The idea of internalization is simple: given a mathematical concept $X$ 
defined solely in terms of sets, functions and commutative diagrams 
involving these, and given some category $K$, one obtains the concept of 
an `$X$ in $K$' by replacing all these sets, functions and commutative
diagrams by corresponding objects, morphisms, and commutative diagrams
in $K$.  

The case we need here is when $X$ is the concept of `category':

\begin{definition} \et 
Let $K$ be a category.  
An {\bf internal category in} $K$, or simply
{\bf category in} $K$, say $C$, consists of: 
\begin{itemize}
\item an object $\Ob(C) \in K$,
\item an object $\Mor(C) \in K$,
\item {\bf source} and {\bf target} morphisms
$s,t \maps \Mor(C) \to \Ob(C)$,
\item an {\bf identity-assigning} morphism
$\id \maps \Ob(C) \to \Mor(C)$,
\item a {\bf composition} morphism
$\circ \maps \Mor(C) {{}_{s}\times_{t}} 
\Mor(C) \to \Mor(C)$
\end{itemize}
satisfying the usual rules of a category expressed in terms of commutative
diagrams.
\end{definition}

\noindent Here $\Mor(C) {{}_{s}\times_{t}} \Mor(C)$ is defined using a
pullback: if $K$ is the category of sets, it is the set of composable
pairs of morphisms in $C$.  Inherent in the definition is the
assumption that this pullback exist, along with the other pullbacks
needed to write the rules of a category as commutative diagrams.

We can similarly define a {\bf functor in} $K$ and a {\bf natural
transformation in} $K$; details can be found in Borceux's handbook
\cite{Borceux}.  There is a 2-category $K\Cat$ whose objects,
morphisms and 2-morphisms are categories, functors and natural
transformations in $K$.  To study \emph{symmetries} in higher
gauge theory, we need these examples:

\begin{definition} \et
\label{Lie 2-group}
Let $\Lie\Grp$ be the category whose objects are Lie groups and whose 
morphisms are Lie group homomorphisms.   Then the objects, morphisms
and 2-morphisms of $\Lie\Grp\Cat$ are called {\bf Lie 2-groups}, 
{\bf Lie 2-group homomorphisms}, and {\bf Lie 2-group 
2-homomorphisms}, respectively.
\end{definition}

\begin{definition} \et
\label{Lie 2-algebra}
Let $\Lie\Alg$ be the category whose objects are Lie algebras and whose 
morphisms are Lie algebra homomorphisms.   Then the objects, morphisms
and 2-morphisms of $\Lie\Alg\Cat$ are called {\bf Lie 2-algebras}, 
{\bf Lie 2-algebra homomorphisms}, and {\bf Lie 2-algebra 
2-homomorphisms}, respectively.
\end{definition}

\noindent
For the benefit of experts, we should admit that we are only
defining `strict' Lie 2-groups and Lie 2-algebras, where all the usual
laws hold as equations.  We rarely need any other kind in this paper, 
but there are more general Lie 2-groups and Lie 2-algebras
where the usual laws hold only up to isomorphism 
\cite{BaezLauda:2003,BaezCrans:2003}.

We could also consider categories in $\Diff$, the category whose objects 
are finite-dimensional smooth manifolds and whose morphisms are smooth maps.  
Ehresmann \cite{Ehresmann:1959} introduced these in the late 1950's
under the name of `differentiable categories'.  However, these 
are not quite what we want here, for two reasons.  First, unlike $\Lie\Grp$ 
and $\Lie\Alg$, $\Diff$ does not have pullbacks in general.   This means
that when we try to define a category in $\Diff$, the set 
of composable pairs of morphisms is not automatically a smooth manifold.
Second, the space of smooth paths in a smooth manifold is not again a 
smooth manifold.  This is an annoyance when studying connections on
bundles.  

To solve these problems, we want a category of `smooth spaces' that has
pullbacks and includes path spaces.  Various categories of this sort 
have been proposed.  It is unclear which one is best, but we shall 
use the last of several variants proposed by Chen \cite{Chen:1977,Chen:1982}.  
In what follows, we use {\bf convex set} to mean a convex subset of 
$\R^n$, where $n$ is arbitrary (not fixed).  Any convex set inherits
a topology from its inclusion in $\R^n$.  We say a map $f$ 
between convex sets is {\bf smooth} if arbitrarily high derivatives 
of $f$ exist and are continuous, using the usual definition
of derivative as a limit of a quotient.   

\begin{definition} \et
\label{smooth_space} 
A {\bf smooth space} is a set $X$ equipped with, for 
each convex set $C$, a collection of functions 
$\phi \maps C \to X$ called {\bf plots} in $X$, such that:
\begin{enumerate} 
\item 
If $\phi \maps C \to X$ is a plot in $X$, 
and $f \maps C' \to C$ is a smooth map between convex sets, 
then $\phi \circ f$ is a plot in $X$,
\item
If $i_\alpha \maps C_\alpha \to C$ is an open cover of a convex
set $C$ by convex subsets $C_\alpha$,
and $\phi \maps C \to X$ has the property that $\phi \circ i_\alpha$ 
is a plot in $X$ for all $\alpha$, then $\phi$ is a plot in $X$.
\item 
Every map from a point to $X$ is a plot in $X$.
\end{enumerate}
\end{definition}

\begin{definition} \et
\label{smooth_map} 
A {\bf smooth map} from the smooth space $X$ to the smooth
space $Y$ is a map $f \maps X \to Y$ such that for every
plot $\phi$ in $X$, $\phi \circ f$ is a plot in $Y$.
\end{definition}

\noindent
In highbrow lingo, this says that smooth spaces are sheaves
on the category whose objects are convex sets and whose morphisms are
smooth maps, equipped with the Grothendieck topology where a cover
is an open cover in the usual sense.  However, smooth spaces are
not arbitrary sheaves of this sort, but precisely those 
for which two plots with domain $C$ agree whenever they agree
when pulled back along every smooth map from a point to $C$.

Using this, it is straightforward to check that there is a category 
${\rm C}^\infty$ whose objects are smooth spaces and whose morphisms
are smooth maps.  Moreover this category is cartesian closed, and 
it has arbitrary limits and colimits.  It also has other nice properties:

\begin{itemize}
\item  Every finite-dimensional smooth manifold (possibly with
boundary) is a smooth space;
smooth maps between these are precisely those
that are smooth in the usual sense.
\item  Every smooth space can be given the strongest topology in which
all plots are continuous; smooth maps
are then automatically continuous.
\item  Every subset of a smooth space is a smooth space.
\item  We can form a quotient of a smooth space $X$ by any
equivalence relation, and the result is again a smooth space.  
\item  We can define vector fields and differential forms on smooth 
spaces, with many of the usual properties.
\end{itemize}

With the notion of smooth space in hand, we can make the following
definitions:

\begin{definition} \et
\label{smooth 2-space}
Let ${\rm C}^\infty$ be the category whose objects are 
smooth spaces and whose morphisms are smooth maps.   Then the objects, 
morphisms and 2-morphisms of ${\rm C}^\infty\Cat$ are called 
{\bf smooth 2-spaces}, {\bf smooth maps}, and {\bf smooth 2-maps}, 
respectively.
\end{definition}

Writing down the above definitions is quick and easy.  It takes 
longer to understand them and apply them to higher gauge theory.
For this we must unpack them and look at examples.   

To get examples of Lie 2-groups, we can use `Lie crossed modules'.  
A {\bf Lie crossed module} is a quadruple $(G,H,t,\alpha)$ 
where $G$ and $H$ are Lie 
groups, $t \maps H \to G$ is a Lie group homomorphism and 
$\alpha$ is a smooth action of $G$ as automorphisms of $H$ such
that $t$ is equivariant:
\[   t(\alpha(g)(h)) = g \, t(h)\, g^{-1} \]
and satisfies the so-called `Peiffer identity':
\[    \alpha(t(h))(h') = hh'h^{-1} . \]
We obtain a Lie crossed module from a Lie 2-group $\twogroup$ as
follows:
\begin{itemize}
\item $G$ is the Lie group of objects of $\G$, $\Ob(\twogroup)$,
\item $H$ is the subgroup of $\Mor(\twogroup)$ consisting
of morphisms with source equal to $1 \in G$:
\item $t \maps H \to G$ is the homomorphism sending each 
morphism in $H$ to its target,
\item $\alpha$ is the action of $G$ as automorphisms of $H$ 
defined using conjugation in $\Mor(\twogroup)$ as follows: 
$\alpha(g) h = 1_g h {1_g}^{-1}$. 
\end{itemize}
Conversely, we can reconstruct any Lie 2-group from its Lie 
crossed module.  In fact, the 2-category of Lie 2-groups is 
biequivalent to that of Lie crossed modules \cite{BaezLauda:2003}.
This gives various examples:

\begin{example} \et
\label{crossed.1} 
{\rm 
Given any abelian group $H$, there is a Lie crossed module
where $G$ is the trivial group and $t, \alpha$ are trivial.  
This gives a Lie 2-group $\twogroup$ with one object
and $H$ as the group of morphisms.  Lie 2-groups of this
sort are important in the theory of \emph{abelian} gerbes.
}\end{example}

\begin{example} \et
\label{crossed.2} 
{\rm
Given any Lie group $H$, there is a Lie
crossed module with $G = \Aut(H)$, $t \maps H \to G$ the homomorphism
assigning to each element of $H$ the corresponding inner automorphism,
and the obvious action of $G$ as automorphisms of $H$.  We call the
corresponding Lie 2-group the {\bf automorphism 2-group} of $H$, and
denote it by $\AUT(H)$.   This sort of 2-group is important in the theory 
of \emph{nonabelian} gerbes.  

We use the term `automorphism 2-group' because $\AUT(H)$ really is the
2-group of symmetries of $H$.  Lie groups form a 2-category, any object
in a 2-category has a 2-group of symmetries, and the 2-group of 
symmetries of $H$ is naturally a Lie 2-group, which is none other than
$\AUT(H)$.  See \cite{BaezLauda:2003} for details.
}\end{example}

\begin{example} \et 
\label{crossed.3} 
{\rm 
Suppose that $1 \to A \hookrightarrow H \stackto{t} G \to 1$ is a central
extension of the Lie group $G$ by the Lie group $H$.  Then there is
a Lie crossed module with this choice of $t \maps H \to G$.  To construct
$\alpha$ we pick any section $s$, that is, any function
$s \maps G \to H$ with $t(s(g)) = g$, and define 
\[      \alpha(g) h = s(g) h s(g)^{-1}  .\]
Since $A$ lies in the center of $H$, $\alpha$ independent of the choice of
$s$.  We do not need a global smooth section $s$ to show $\alpha(g)$ depends 
smoothly on $g$; it suffices that there exist a local smooth section in a
neighborhood of each $g \in G$.

It is easy to generalize this idea to infinite-dimensional cases
if we work not with Lie groups
but {\bf smooth groups}: that is, groups in the category of smooth spaces.  
The basic theory of smooth groups, smooth 2-groups and smooth crossed 
modules works just like the finite-dimensional case, but with the category 
of smooth spaces replacing $\Diff$.  In particular, every smooth group 
$G$ has a Lie algebra $\g$.  

Given a connected and simply-connected compact simple Lie 
group $G$, the loop group $\Omega G$ is a smooth group. 
For each `level' $k \in \Z$, this group has a central extension
\[           1 \to \U(1) \hookrightarrow 
\widehat{\Omega_k G} \stackto{t} \Omega G \to 1   \]
as explained by Pressley and Segal \cite{PressleySegal:1986}.  The above 
diagram lives in the category of smooth groups, and there 
exist local smooth sections for $t \maps \widehat{\Omega_k G} \to \Omega G$, 
so we obtain a smooth crossed module $(\Omega G, \widehat{\Omega_k G}, 
t,\alpha)$ with $\alpha$ given as above.  This in turn gives an
smooth 2-group which we call the \textbf{level-\textit{k} 
loop 2-group} of $G$, ${\cal L}_k G$.

It has recently been shown \cite{BCSS}
that ${\cal L}_k G$ fits into an exact sequence of smooth 2-groups:
\[           1 \to {\cal L}_k G \hookrightarrow 
{\cal P}_k G \stackto{} G \to 1   \]
where the middle term, the \textbf{level-\textit{k} path 2-group} 
of $G$, has very interesting properties.  In particular, when $k = \pm 1$,
the geometric realization of the nerve of ${\cal P}_k G$ is a topological 
group that can also be obtained by killing the 3rd homotopy group of $G$.
When $G = {\rm Spin}(n)$, this topological group goes by the name of
${\rm String}(n)$, since it plays a role in defining spinors on loop
space \cite{Witten:1988}.  The group ${\rm String}(n)$ also shows up in 
Stolz and Teichner's work on elliptic cohomology, which involves a 
notion of parallel transport over surfaces \cite{StolzTeichner}.  
So, we expect that ${\cal P}_k G$ will be an especially interesting 
structure 2-group for applications of 2-bundles to string theory.

To define the holonomy of a connection, we need smooth groups with 
an extra property: namely, that for every smooth function 
$f \maps [0,1] \to \g$ there is a unique smooth function 
$g \maps [0,1] \to G$ solving the differential equation
\[       {d\over dt}g(t) = f(t)g(t)  \]
with $g(0) = 1$.  We call such smooth groups {\bf exponentiable}.
Similarly, we call a smooth 2-group $\twogroup$ {\bf exponentiable} 
if its crossed module $(G,H,t,\alpha)$ has both $G$ and $H$ 
exponentiable.  In particular, every Lie group and thus every
Lie 2-group is exponentiable.  The smooth groups $\Omega G$
and $\widehat{\Omega_k G}$ are also exponentiable, as are the 
2-groups ${\cal L}_k G$ and ${\cal P}_k G$.  So, for the convenience
of stating theorems in a simple way, {\it we henceforth implicitly
assume all smooth groups and 2-groups under discussion are exponentiable}.
We only really need this in Theorems \ref{central proposition concerning
connections} and \ref{central proposition concerning 2-connections}.
}\end{example}

Finally, here are some easy examples of smooth 2-spaces:

\begin{example} \et 
\label{2space.1}
{\rm 
Any smooth space can be seen as a 
smooth 2-space with only identity morphisms.
}
\end{example}

\begin{example} \et 
\label{2space.2} 
{\rm 
Any smooth group (for example a Lie group) can be seen as
a smooth 2-space with only one object.
}
\end{example}

\begin{example} \et \label{2space.3}
{\rm 
Given a smooth space $M$, there is a smooth 2-space $\P_1(M)$,
the {\bf path groupoid of $M$}, such that:
\begin{itemize}
\item the objects of $\P_1(M)$ are points of $M$,
\item the morphisms of $\P_1(M)$ are thin homotopy classes of smooth paths
$\gamma \maps [0,1] \to M$ such that $\gamma(s)$ is constant near $s = 0$
and $s = 1$.
\end{itemize}
Here a {\bf thin homotopy} between smooth paths $\gamma_1, \gamma_2
\maps [0,1] \to M$ is a smooth map $H \maps [0,1]^2 \to M$ such that:
\begin{itemize}
\item 
$H(s,0) = \gamma_1(s)$ and 
$H(s,1) = \gamma_2(s)$, 
\item
$H(s,t)$ is independent of $t$ near $t = 0$ and near $t = 1$, 
\item
$H(s,t)$ is constant near $s = 0$ and near $s = 1$,
\item
the rank of the differential $dH(s,t)$ is $< 2$ for all $s,t \in
[0,1]$. 
\end{itemize}
The last condition is what makes the homotopy `thin': it guarantees
that the homotopy sweeps out a surface of vanishing area.  

To see how $\P_1(M)$ becomes a 2-space, first note that
the space of smooth maps $\gamma \maps [0,1] \to M$ becomes a
smooth space in a natural way, as does the subspace satisfying
the constancy conditions near $t = 0,1$, and finally the quotient
of this subspace by the thin homotopy relation.  
This guarantees that $\Mor(\P_1(M))$ is a smooth space.  Clearly
$\Ob(\P_1(M)) = M$ is a smooth space as well.  One can check that
$\P_1(M)$ becomes a smooth 2-space with usual composition of paths
giving the composition of morphisms.

In fact, $\P_1(M)$ is not just a smooth 2-space: it is also a groupoid.
The inverse of $[\gamma]$ is just $[\overline{\gamma}]$, 
where $\overline{\gamma}$ is obtained by reversing the orientation 
of the path $\gamma$:
\[        \overline{\gamma}(s) = \gamma(1 - s)  .\]
Moreover, the map sending any morphism to its inverse is 
smooth.  Thus $\P_1(M)$ is a {\bf smooth groupoid}: a 2-space
where every morphism is invertible and the map sending every morphism
to its inverse is smooth.  
}\end{example}  

\section{2-Bundles}
\label{section: 2-Bundles}

In differential geometry an
ordinary bundle consists of two smooth spaces, the {\bf total space}
$E$ and the {\bf base space} $B$, together with a {\bf projection map} 
\[
  E \stackto{p} B
  \,.
\]
To categorify the theory of bundles, we start by replacing smooth spaces by 
smooth 2-spaces:

\begin{definition} \et
  \label{2-bundle}
  A {\bf 2-bundle} consists of
  \begin{itemize}
    \item
      a smooth 2-space $E$ (the {\bf total 2-space}),
    \item
      a smooth 2-space $B$ (the {\bf base 2-space}),
    \item
      a smooth map $p \maps E \to B$ 
      (the {\bf projection}).
  \end{itemize}
\end{definition} 

In gauge theory we are interested in \emph{locally trivial} 2-bundles. 
Ordinarily, a locally trivial
bundle with fiber $F$ is a bundle $E \stackto{p} B$ together with an open
cover $U_i$ of $B$, such that the restriction of $E$ to any of the $U_i$ 
is equipped with an isomorphism to the trivial bundle $U_i \times F \to U_i$. 
To categorify this, we would need to define a `2-cover' of the base 2-space 
$B$.  This is actually a rather tricky issue, since forming the `union' of 
2-spaces requires knowing how to compose a morphism in one 2-space
with a morphism in another.  While this issue can be addressed, we prefer 
to avoid it here by assuming that $B$ is just an ordinary smooth space,
regarded as a smooth 2-space with only identity morphisms.

We can now state the definition of a
locally trivial 2-bundle.  First note that we can restrict a
2-bundle $E \stackto{p} B$ to any subspace $U \subseteq B$
to obtain a 2-bundle which we denote by $E|_U \stackto{p} U$.  Then:

\begin{definition} \et
\label{2-bundle with local trivialization}
Given a smooth 2-space $F$, we define a {\bf locally trivial 2-bundle 
with fiber $F$} to be a 2-bundle $E \stackto{p} B$ and an open cover 
$\{U_i\}$ of the base space $B$ 
equipped with equivalences 
  \[
    E|_{\twoU_i} \stackto{t_i} \twoU_i \times \twoF
  \]
called {\bf local trivializations} such that these diagrams: 
\[ 
\xymatrix @!0
{ E|_{\twoU_i}
 \ar [dddrr]_{p}
  \ar[rrrr]^{t_i} 
  & &  & &
  U_i \times F
  \ar[dddll]^{}
  \\ \\ &  \\ & &
   U_i 
 }
\]
commute for all $i \in I$.
\end{definition}
Readers wise in the ways of categorification \cite{BaezDolan:1998} may 
ask why we did not require that these diagrams commute {\it up to 
natural isomorphism}.  The reason is that $U_i$, as an ordinary space, has
only identity morphisms when we regard it as a 2-space.  Thus, for
this diagram to commute up to natural isomorphism, it must commute
`on the nose'. 

Readers less wise in the ways of categorification may find the
above definition painfully abstract.  So, let us translate it
into data that specify how to build a locally trivial 2-bundle
from trivial ones over the patches $U_i$.  For this, 
we need to extract \emph{transition functions} from the local 
trivializations.

So, suppose $E \stackto{p} B$ is a locally trivial 2-bundle with 
fiber $F$.  This means that $B$ is equipped with an 
open cover $U$ and for each open set $U_i$ in the cover
we have a local trivialization 
\[   t_i \maps  E|_{U_i} \to \twoU_i \times \twoF \]
which is an equivalence.  This means
that $t_i$ is equipped with a specified map
\[   \bar t_i  \maps \twoU_i \times \twoF \to E|_{\twoU_i} \]
together with invertible 2-maps
\[  
\begin{array}{ccl}
   \tau_i \maps \bar t_i t_i &\To& 1  \\
   \bar \tau_i \maps t_i \bar t_i &\To& 1   
\end{array}
\] 
In particular, this means that $\bar t_i$ is also an equivalence.

Now consider a double intersection $U_{ij} = U_i \cap U_j$.
The composite of equivalences is again an equivalence, so
we get an {\bf autoequivalence}
\[
t_j \bar t_i
\maps U_{ij}\times \twoF \to U_{ij}\times \twoF
\]
that is, an equivalence from this 2-space to itself.
By the commutative
diagram in \refdef{2-bundle with local trivialization}, this
autoequivalence 
must act trivially on the $U_{ij}$ factor, so 
\[
  t_j \bar t_i(x,f) = (x,f g_{ij}(x)) 
\]
for some smooth function $g_{ij}$ from $U_{ij}$ to the 
smooth space of autoequivalences of the fiber $F$.  
Note that we write these autoequivalences as acting on 
$F$ from the right, as customary in the theory of bundles.
We call the functions $g_{ij}$ {\bf transition functions}, since they
are just categorified versions of the usual transition functions 
for locally trivial bundles.  

In fact, for any smooth 2-space $F$ there is a smooth 2-space $\AUT(F)$ 
whose objects are autoequivalences of $F$ and whose morphisms are 
invertible 2-maps between these.  The transition functions are maps
\[     g_{ij} \maps U_{ij} \to \Ob(\AUT(F))  .\]
The 2-space $\AUT(F)$ is a kind of 2-group, with composition of
autoequivalences giving the product.  However, is not the sort of
2-group we have been considering here, because it does not have
`strict inverses': the group laws involving inverses do not hold as
equations, but only up to specified isomorphisms that satisfy
coherence laws of their own.  So, $\AUT(F)$ is a `coherent' smooth
2-group in the sense of Baez and Lauda \cite{BaezLauda:2003}.

Next, consider a triple intersection $U_{ijk} = U_i \cap U_j \cap U_k$.
In an ordinary locally trivial bundle the transition functions satisfy 
the equation $g_{ij} g_{jk} = g_{ik}$, but in a locally trivial 
2-bundle this holds {\it only up to isomorphism}.  In other
words, there is a smooth map 
\[    h_{ijk} \maps U_{ijk} \to \Mor(\AUT(F))  \]
such that for any $x \in U_{ijk}$, 
\[    h_{ijk}(x) \maps g_{ij}(x) g_{jk}(x) \stackto{\sim} g_{ik}(x) . \]
To see this, note that there is an invertible 2-map
\[ t_k \tau_j \bar t_i \, \maps \,
t_k \bar t_j t_j \bar t_i \To t_k \bar t_i  \]
defined by horizontally composing $\tau_j$ with $t_k$ on the left
and $\bar t_i$ on the right.   Since
\[   t_k \bar t_j t_j \bar t_i(x,f) = (x,fg_{ij}(x)g_{jk}(x)) \]
while
\[     t_k \bar t_i(x,f) = (x, f g_{ik}(x))  \]
we have
\[ t_k \tau_j \bar t_i (x,f) \maps 
(x, fg_{ij}(x)g_{jk}(x)) \to (x,fg_{ik}(x)).\]
Since this morphism must be the identity on the first factor,
we have
\[ t_k \tau_j \bar t_i (x,f) = (1_x, f h_{ijk}(x)) \]
where $h_{ijk}(x) \maps g_{ij}(x) g_{jk}(x) \to g_{ik}(x)$ depends
smoothly on $x$. 

Similarly, in a locally trivial bundle we have
$g_{ii} = 1$, 
but in a locally trivial 2-bundle there is a smooth map
\[    k_i \maps U_i \to \Mor(\AUT(F))  \]
such that for any $x \in U_i$,
\[    k_i(x) \maps g_{ii}(x) \to 1  .\]
To see this, recall that there is an invertible 2-map
\[     \bar\tau_i \maps t_i \bar t_i \To 1 . \]
Since
\[   t_i \bar t_i(x,f) = (x, f g_{ii}(x)) \]
we have
\[ \bar \tau_i (x,f) \maps (x, f g_{ii}(x)) 
\to (x,f) , \]
and since this morphism must be the identity on the first factor,
we have
\[ \bar \tau_i (x,f) = (1_x, f k_i(x)) \]
where $k_i(x) \maps g_{ii}(x) \to 1$ depends smoothly on $x$. 

In short, the transition functions $g_{ij}$ for a locally trivial 2-bundle
satisfy the usual cocycle conditions up to specified isomorphisms
$h_{ijk}$ and $k_i$, which we call {\bf higher transition functions}.
These, in turn, satisfy some cocycle conditions of their own: 

\begin{theorem} \et 
\label{2-bundle cocycle conditions} 
Suppose $E \stackto{p} B$ is a locally trivial 2-bundle, 
and define the transition functions 
$g_{ij}, h_{ijk}$, and $k_i$ as above.  Then:

\begin{itemize}
\item
$h$ makes this diagram,
called the {\bf associative law},
commute for any $x \in U_{ijkl}$:
\[
\xy
 (0,20)*+{g_{ij}(x)\, g_{kl}(x) \, g_{lm}(x) }="T";
 (-20,0)*+{g_{il}(x)\, g_{lm}(x)}="L";
 (20,0)*+{g_{ij}(x)\, g_{jm}(x)}="R";
 (0,-20)*+{g_{im}(x)}="B";
 (-7,0)*{}="ML";
 (7,0)*{}="MR";
     {\ar^{g_{ij}(x) \, h_{jkl}(x)} "T";"R"};
     {\ar_{h_{ijk}(x) \, g_{kl}(x)} "T";"L"};
     {\ar^{h_{ijl}(x)} "R";"B"};
     {\ar_{h_{ikl}(x)} "L";"B"};
\endxy
\]

\item 
$k$ makes these diagrams, 
called the {\bf left and right unit laws}, 
commute for any $x \in U_{ij}$:
\[
\xy
  (-36,5)*+{g_{ii}(x) \, g_{ij}(x) }="L";
  (0,5)*+{1 \, g_{ij}(x)}="M";
  (0,-18)*+{g_{ij}(x)}="B";
  (-5,2)*{}="TL";
  (5,2)*{}="TR";
  (10,-4)*{}="BR";
  (-10,-4)*{}="BL";
     {\ar^{k_i(x) \, g_{ij}(x)} "L";"M"};
     {\ar^{=} "M";"B"};
     {\ar_{h_{iij}(x)} "L";"B"};
\endxy
\qquad
\xy
  (0,5)*+{g_{ij}(x) \, 1}="M";
  (36,5)*+{g_{ij}(x) \, g_{jj}(x)}="R";
  (0,-18)*+{g_{ij}(x)}="B";
  (-5,2)*{}="TL";
  (5,2)*{}="TR";
  (10,-4)*{}="BR";
  (-10,-4)*{}="BL";
     {\ar_{g_{ij}(x) \, k_j(x)} "R";"M"};
     {\ar_{=} "M";"B"};
     {\ar^{h_{ijj}(x)} "R";"B"};
\endxy
\]
\end{itemize}
\end{theorem}

\Proof  
Checking that these diagrams commute is a straightforward 
computation using the definitions of $g,h$, and $k$ in 
terms of $t,\bar t$, $\tau$ and $\bar \tau$.  
\endofproof

The associative law and unit laws are analogous to those which hold 
in a monoid.  They also have simplicial interpretations.  
In a locally trivial bundle, the transition functions
give a commuting triangle for any triple intersection:
\begin{center}
  \begin{picture}(100,85)
  \includegraphics{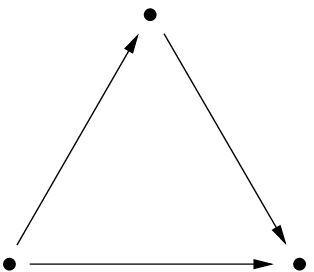}
  \put(-80,40){$g_{ij}$}
  \put(-24,40){$g_{jk}$}
  \put(-50,-7){$g_{ik}$}
  \end{picture}
\end{center}
\vskip 2em 
In a locally trivial 2-bundle, such triangles commute only
up to isomorphism:
\begin{center}
 \begin{picture}(100,85)
  \includegraphics{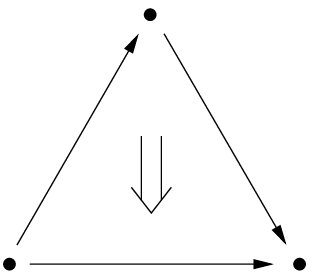}
  \put(-120,0){}
  \put(-60,80){}
  \put(2,0){}
  \put(-80,40){$g_{ij}$}
  \put(-24,40){$g_{jk}$}
  \put(-50,-7){$g_{ik}$}
  \put(-70,24){$h_{ijk}$}
  \end{picture}
  \vskip 1em
\end{center}
However, the associative law says that for each quadruple 
intersection, this tetrahedron commutes:
\begin{eqnarray*}
\begin{picture}(140,160)
  \includegraphics{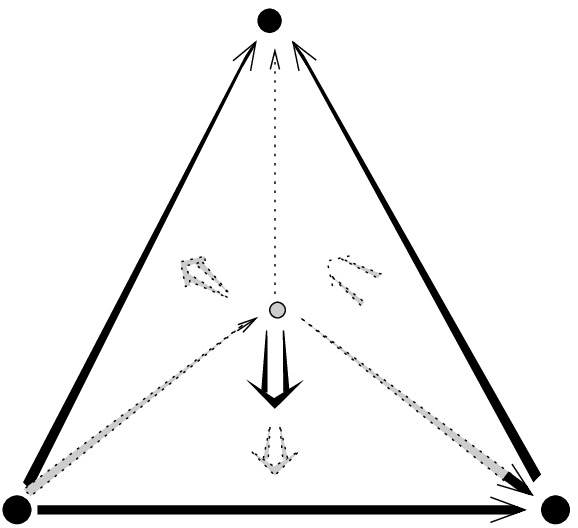}
  \put(-127,46){${}_{g_{ij}}$}
  \put(-53,46){${}_{g_{jk}}$}
  \put(-90,-3){$g_{ik}$}
  \put(-137,80){$g_{il}$}
  \put(-41,80){$g_{kl}$}
  \put(-82,105){${}_{g_{jl}}$}

  \put(-80,22){$h_{ijk}$}
  \put(-82,47){$h_{ikl}$}
  \put(-68,85){${}_{h_{jkl}}$}
  \put(-112,85){${}_{h_{ijl}}$}
\end{picture}
\end{eqnarray*}

\vskip 1em \noindent

\noindent
We can also visualize the left and right unit laws simplicially,
but they involve degenerate tetrahedra:

\begin{center}
\begin{picture}(230,180)
\includegraphics{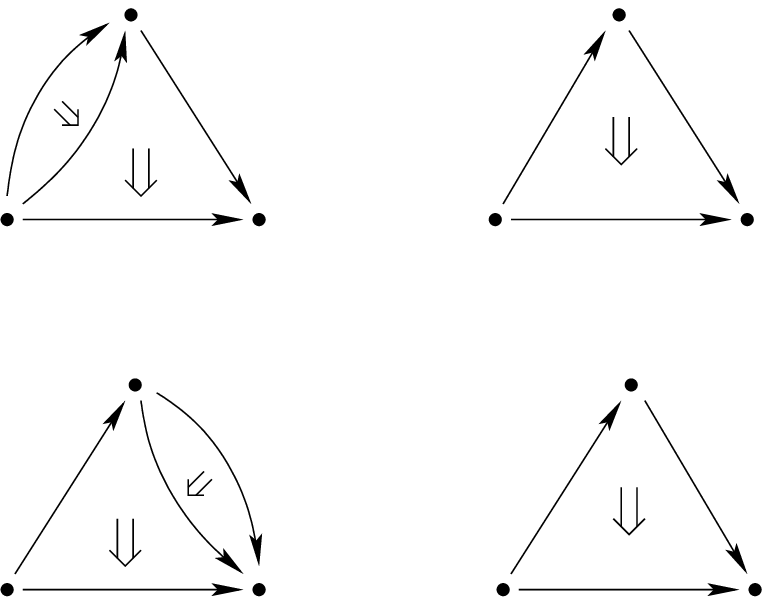}
\put(0,0){}
\put(-45,-6){$g_{ij}$}
\put(-190,-6){$g_{ij}$}
\put(-45,100){$g_{ij}$}
\put(-190,100){$g_{ij}$}
\put(-20,145){$g_{ij}$}
\put(-161,145){$g_{ij}$}
\put(-71,35){$g_{ij}$}
\put(-215,35){$g_{ij}$}
\put(-75,145){$g_{ii}$}
\put(-224,145){$g_{ii}$}
\put(-200,144){$k_i$}
\put(-193,129){$1$}
\put(-176,122){$1$}
\put(-117,130){$=$}
\put(-55,120){$h_{iij}$}
\put(-117,20){$=$}
\put(-16,35){$g_{jj}$}
\put(-149,35){$g_{jj}$}
\put(-172,39){$k_j$}
\put(-175,20){$1$}
\put(-194,15){$1$}
\put(-52,13){$h_{ijj}$}
\end{picture}
\end{center}

We are now almost in a position to define $\twogroup$-2-bundles 
for any smooth 2-group $\twogroup$; we only need to understand
how a 2-group can `act' on a 2-space.  For simplicity we only
consider the case of a strict action:

\begin{definition} \et 
A {\bf (strict) action} of a smooth 2-group 
$\twogroup$ on a smooth 2-space $F$ is a smooth homomorphism
\[    \alpha \maps \twogroup \to \AUT(F) , \]
that is, a smooth map that preserves products and inverses.
\end{definition}

\noindent
Note in particular that every smooth 2-group has an action on itself
via right multiplication.

\begin{definition} \et
For any smooth 2-group $\twogroup$, we say a locally trivial
2-bundle $E \to M$ has $\twogroup$ as its {\bf structure 
2-group} when the transition functions $g_{ij}$, $h_{ijk}$, and 
$k_i$ factor through an action $\twogroup \to \AUT(F)$.  In this case we 
also say $P$ is a {\bf $\twogroup$-2-bundle}.
If furthermore $F = \twogroup$ and $\twogroup$ acts on $F$ by right
multiplication, we say $P$ is a {\bf principal $\twogroup$-2-bundle}.  
\end{definition}

For a principal $\twogroup$-2-bundle we can think of the transition 
functions as taking values in the groups $\Ob(\twogroup)$ and 
$\Mor(\twogroup)$.  The reader familiar with 
gerbes will note that these functions, satisfying the equations they do, 
reduce to the usual sort of cocycle defining an 
\emph{abelian gerbe} when $k_i = 1$ and $\twogroup$ has
the special form described in Example \ref{crossed.1}.  Similarly, 
they reduce to a cocycle defining a \emph{nonabelian gerbe}
when $k_i = 1$ and $\twogroup$ has the form described in 
Example \ref{crossed.2}.  Thus there is a close relation between
principal 2-bundles and gerbes, much like that between principal
bundles and sheaves of groups.

The equation $k_i = 1$ arises because gerbes 
are often defined using {\v C}ech cocycles that are antisymmetric 
in the indices $i,j,k,\dots,$ in the sense that group-valued functions 
go to their inverses upon an odd permutation of these indices.  Thus in 
this context $g_{ii} = 1$, and one implicitly assumes $k_i = 1$.    
In fact, Bartels \cite{Bartels:2004} has shown that every 
$\twogroup$-2-bundle is equivalent to one with $k_i = 1$.  To 
state this result, he first needed to define a 2-category of 
$\twogroup$-2-bundles.  This 2-category is equivalent to the 
2-category of abelian or nonabelian gerbes when $\twogroup$ has 
one of the two special forms mentioned above.

\section{2-Connections}
\label{section: 2-Connections}

For a trivial bundle, 
the holonomy of a connection assigns elements of the structure
group to paths in space.  Similarly, a 2-connection
assigns objects and morphisms of the structure 2-group to 
paths and surfaces in space.  To make this precise we need 
the notion of a `path 2-groupoid'.  

We described the path groupoid of a smooth space $M$
in Example \ref{2space.3}.  This has points of $M$ as
objects:
\[    \bullet \, {\textstyle{\small{\textit{x}}}}  \]
and thin homotopy classes of paths in $M$ as morphisms:
\[
\xymatrix{
 {\textstyle{\small{\textit{x}}}}  \,
   \bullet \ar@/^1pc/[rr]^{\gamma}
&& \bullet \, {\textstyle{\small{\textit{y}}}}   
}
\]
The path 2-groupoid also has 2-morphisms, which are
thin homotopy classes of 2-dimensional surfaces like this: 
\[
\xymatrix{
 {\textstyle{\small{\textit{x}}}} \, 
   \bullet \ar@/^1pc/[rr]^{\gamma_1}_{}="0"
           \ar@/_1pc/[rr]_{\gamma_2}^{}="1"
&& \bullet \, {\textstyle{\small{\textit{y}}}}   
\ar@{=>}"0";"1"^{\Sigma}
}
\]
We call these `bigons':

\begin{definition} \et 
\label{parametrized bigon}
  Given a smooth space $M$, a {\bf parametrized bigon} in $M$
  is a smooth map
  \[
    \Sigma \maps [0,1]^2 \to M
  \]
  which is constant near $s = 0$, constant near $s = 1$, 
  independent of $t$ near $t = 0$, and independent of $t$ near $t = 1$.
  We call $\Sigma\of{\cdot,0}$ the {\bf source} of the 
  parametrized bigon $\Sigma$, and $\Sigma\of{\cdot,1}$ the 
  {\bf target}.  If $\Sigma$ is a parametrized bigon with
  source $\gamma_1$ and target $\gamma_2$, we write $\Sigma \maps
  \gamma_1 \to \gamma_2$.   
\end{definition}

\begin{definition} \et
\label{bigon}
Suppose $\Sigma \maps \gamma_1 \to \gamma_2$ and
$\Sigma' \maps \gamma'_1 \to \gamma'_2$ are parametrized
bigons in a smooth space $M$.  A {\bf thin homotopy} between
$\Sigma$ and $\Sigma'$ is a smooth map 
\[
  H \maps [0,1]^3 \to M
\]
with the following properties:
\begin{itemize}

\item 
$H\of{s,t,0} = \Sigma\of{s,t}$ and $H\of{s,t,1} = \Sigma'\of{s,t}$,

\item 
$H\of{s,t,u}$ is independent of $u$ near $u = 0$ and near $u = 1$,

\item 
For some thin homotopy $F_1$ from $\gamma_1$ to $\gamma_1'$,
$H\of{s,t,u} = F_1\of{s,u}$ for $t$ near $0$, and
for some thin homotopy $F_2$ from $\gamma_2$ to $\gamma_2'$, 
$H\of{s,t,u} = F_2\of{s,u}$ for $t$ near $1$,

\item 
$H(s,t,u)$ is constant for $s = 0$ and near $s = 1$,

\item 
$H$ does not sweep out any volume: the rank of the
differential $dH(s,t,u)$ is $< 3$ for all $s,t,u \in [0,1]$.

\end{itemize}

\noindent
We say two parametrized bigons $\Sigma,\Sigma'$
lie in the same {\bf thin homotopy class}
if the pair $(\Sigma, \Sigma')$ lies in the closure of the 
thin homotopy equivalence relation.  
A {\bf bigon} is a thin homotopy class $[\Sigma]$ of parametrized bigons.
\end{definition}

\begin{definition} \et
\label{2-groupoid of bigons}
  The {\bf path 2-groupoid} 
  $\P_2\of{M}$
  of a smooth space $M$ is the 2-category in which:
  \begin{itemize}
    \item 
      objects are points $x \in M$:
      \[   \bullet \, {\textstyle{\small{\textit{x}}}}  
      \] 
    \item 
      morphisms are thin homotopy classes of paths $\gamma$ in
      $M$ that are constant near $s = 0$ and $s = 1$: 
     \[
      \xymatrix{
       x \ar@/^1pc/[rr]^{[\gamma]} 
       && y 
      }
      \]
    \item
      2-morphisms are bigons in $M$
 \[
\xymatrix{
   x \ar@/^1pc/[rr]^{[\gamma_1]}_{}="0"
           \ar@/_1pc/[rr]_{[\gamma_2]}_{}="1"
           \ar@{=>}"0";"1"^{{}_{[\Sigma]}}
&& y
}
\]
  \end{itemize}
and whose composition operations are defined as:
\begin{itemize}

\item   \hskip 3em
$
\xymatrix{
   x \ar@/^1pc/[rr]^{[\gamma_1]}
&& y \ar@/^1pc/[rr]^{[\gamma_2]}
&& z
}
  =
\xymatrix{
   x \ar@/^1pc/[rr]^{[\gamma_1 \circ \gamma_2]}
&& z 
}
$ 
\hfill\break
\vskip 1em \noindent
where
 \[ 
    \left(\gamma_1\circ \gamma_2\right)\of{s}
    :=
    \left\lbrace
      \begin{array}{cl}
         \gamma_1\of{2s} & \mbox{for $0 \leq s \leq 1/2$}\\
         \gamma_2\of{2s-1} & \mbox{for $1/2\leq s \leq 1$}
      \end{array}
      \right.
\]

\item \hskip 3em
$
\xymatrix{
   x \ar@/^2pc/[rr]^{[\gamma_1]}_{}="0"
           \ar[rr]^<<<<<<{[\gamma_2]}_{}="1"
           \ar@{=>}"0";"1"^{[\Sigma_1]}
           \ar@/_2pc/[rr]_{[\gamma_3]}_{}="2"
           \ar@{=>}"1";"2"^{[\Sigma_2]}
&& y
}
=
\xymatrix{
   x \ar@/^1pc/[rr]^{[\gamma_1]}_{}="0"
           \ar@/_1pc/[rr]_{[\gamma_3]}_{}="1"
           \ar@{=>}"0";"1"^{{}_{[\Sigma_1 \Sigma_2]}}
&& y
}
$
\hfill\break
\vskip 1em \noindent
where
\[
    \left(\Sigma_1\Sigma_2\right)\of{s,t}
    :=
    \left\lbrace
      \begin{array}{cl}
         \Sigma_1\of{s,2t} & \mbox{for $0 \leq t \leq 1/2$}\\
         \Sigma_2\of{s,2t-1} & \mbox{for $1/2\leq t \leq 1$}
      \end{array}
     \right.
\]
    
\item  \hskip 3em
$
\xymatrix{
   x \ar@/^1pc/[rr]^{[\gamma_1]}_{}="0"
           \ar@/_1pc/[rr]_{\gamma_1^\prime}_{}="1"
           \ar@{=>}"0";"1"^{[\Sigma_1]}
&& y \ar@/^1pc/[rr]^{[\gamma_2]}_{}="2"
           \ar@/_1pc/[rr]_{[\gamma_2^\prime]}_{}="3"
           \ar@{=>}"2";"3"^{[\Sigma_2]}
&& z
}
 =
\xymatrix{
   x \ar@/^1pc/[rr]^{[\gamma_1 \circ \gamma_2]}_{}="0"
           \ar@/_1pc/[rr]_{[\gamma_1^\prime\circ\gamma_2^\prime]}_{}="1"
           \ar@{=>}"0";"1"^{{}_{[\Sigma_1 \circ\Sigma_2]}}
&& z 
}
$
\hfill\break
\vskip 1em \noindent
where
\[
    \left(\Sigma_1\circ\Sigma_2\right)\of{s,t}
    :=
    \left\lbrace
      \begin{array}{cl}
         \Sigma_1\of{2s,t} & \mbox{for $0 \leq s \leq 1/2$}\\
         \Sigma_2\of{2s-1,t} & \mbox{for $1/2\leq s \leq 1$}
      \end{array}
    \right.  
\]

  \end{itemize}
\end{definition}

\noindent
One can check that these operations are well-defined, where
for vertical composition we must choose suitable representatives 
of the bigons being composed.  One can also check that
$\P_2(M)$ is indeed a 2-category.  Furthermore, 
the objects, morphisms and 2-morphisms in $\P_2(M)$ all form 
smooth spaces, by an elaboration of the ideas
in Example \ref{2space.3}, and all the 2-category operations are then 
smooth maps.   We thus say $\P_2(M)$ is a {\bf smooth 2-category}:
that is, a 2-category in ${\rm C}^\infty$.  Indeed, the usual definitions 
\cite{KellyStreet} of 2-category, 2-functor, pseudonatural transformation, 
and modification can all be internalized in ${\rm C}^\infty$, and we 
use these `smooth' notions in what follows.
Furthermore, both morphisms and 2-morphisms in $\P_2(M)$ have strict 
inverses, and the operations of taking inverses are smooth, so we say
$\P_2(M)$ is a {\bf smooth 2-groupoid}. 

We obtain the notion of `2-connection' by categorifying
the concept of connection.  The following result 
suggests a strategy for doing this:

\begin{theorem} \et 
\label{central proposition concerning connections}
   For any smooth group $G$ and smooth space $B$, 
   suppose $E \to B$ is a principal $G$-bundle equipped with local
   trivializations over open sets $\{U_i\}_{i \in I}$ covering $B$.
   Let $g_{ij}$ be the transition functions.
   Then there is a one-to-one correspondence between connections on
   $E$ and data of the following sort:
\begin{itemize}
\item 
       for each $i \in I$ a smooth map between smooth 2-spaces:
\[
    \hol_i \maps \P_1\of{U_i} \to G
\]
 called the {\bf local holonomy functor},
 from the path groupoid of $U_i$ 
 to the group $G$ regarded as a smooth 2-space with a single
 object $\bullet$,
\end{itemize}
such that:
\begin{itemize} 
\item
  for each $i,j\in I$, the transition function $g_{ij}$
  defines a smooth natural isomorphism:
  \[
    \hol_i|_{U_{ij}} \stackto{g_{ij}} \hol_j|_{U_{ij}}
  \]
  called the {\bf transition natural isomorphism}.   
  In other words, this diagram commutes:
\begin{center}
\begin{picture}(100,100)
  \includegraphics{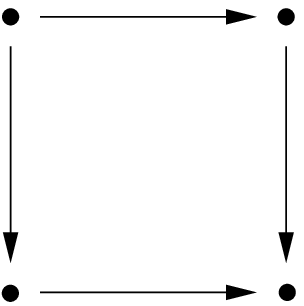}
  \put(-60,90){$g_{ij}\of{x}$}
  \put(-60,-5){$g_{ij}\of{y}$}
  \put(-120,42){$\hol_i(\gamma)$}
  \put(0,42){$\hol_j(\gamma)$}
\end{picture}
\vskip 1em
\end{center}
\noindent
for any path $\gamma \maps x \to y$ in $U_{ij}$.
\end{itemize}
\end{theorem}

\Proof 
See Baez and Schreiber \cite{BaezSchreiber:2004}.
\endofproof

\noindent
In addition, it is worth noting that whenever we have a connection, 
for each $i,j,k \in I$ this triangle commutes:

\begin{center}
  \begin{picture}(70,75)
  \includegraphics{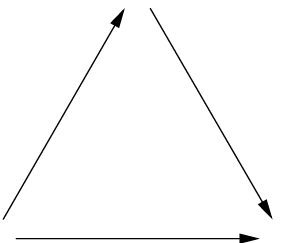}
  \put(-78,40){$g_{ij}$}
  \put(-18,40){$g_{jk}$}
  \put(-50,-7){$g_{ik}$}
  \put(-97,-5){$\hol_i$}
  \put(-50,75){$\hol_j$}
  \put(0,-5){$\hol_k$}
  \end{picture}
\end{center}
\vskip 1em 

The idea behind the above result is that:

  \begin{itemize}
    \item
      The local holonomy functors $\hol_i$ are specified by
      1-forms 
      \[
        A_i \in \Omega^1\of{U_i,\g}\,.
      \]
    \item
      The transition natural isomorphisms $g_{ij}$ 
      are specified by smooth functions
      \[
        g_{ij} \maps U_{ij} \to G\,,
      \] 
      satisfying the equation
      \[
        A_i = g_{ij}A_j g_{ij}^{-1} + g_{ij} \extd g_{ij}^{-1}
      \]
      on $U_{ij}$.
   \item
    The commuting triangle for the triple intersection
    $U_{ijk}$ is equivalent to the equation
    \[
     g_{ij}g_{jk} = g_{ik} 
    \]
     on $U_{ijk}$.
  \end{itemize}

Categorifying all this, we make the following definition:

\begin{definition} \et
  \label{def: 2-connection on a 2-bundle}  
   For any smooth 2-group $\twogroup$, 
   suppose that $E \to B$ is a principal $\twogroup$-2-bundle equipped with 
   local trivializations over open sets $\{U_i\}_{i \in I}$ covering $B$, 
   and let the transition functions $g_{ij}$, 
   $h_{ijk}$ and $k_i$ be given as in Theorem
   \ref{2-bundle cocycle conditions}.  Suppose for simplicity that
   $k_i = 1$.  Then a {\bf 2-connection} on $E$
   consists of the following data:
   \begin{itemize}
\item
   for each $i \in I$ a smooth 2-functor
\begin{eqnarray*}
  \begin{array}{rcccc}
    \hol_i &\maps& \P_2\of{U_i} &\to& \twogroup
    \\
    &&
\xymatrix{
   x \ar@/^1pc/[rr]^{\gamma}="0"
           \ar@/_1pc/[rr]_{}_{\eta}="1"
           \ar@{=>}"0";"1"^{\Sigma}
&& y
}
&\mapsto& 
\xymatrix{
   \bullet \ar@/^1pc/[rr]^{\hol_i\of{\gamma}}="0"
           \ar@/_1pc/[rr]_{}_{\hol_i\of{\eta}}="1"
           \ar@{=>}"0";"1"^{\hol_i\of{\Sigma}}
&& \bullet
}
  \end{array}
\end{eqnarray*}
  called the {\bf local holonomy 2-functor},
  from the path 2-groupoid $\P_2\of{U_i}$ 
  to the 2-group $\twogroup$ regarded as a smooth 2-category with a single 
  object $\bullet$, 
\end{itemize}
such that:
\begin{itemize}
\item
  For each $i,j$ a pseudonatural isomorphism:
  \[
   g_{ij} \maps \hol_i|_{\P(U_i \cap U_j)}  \to
   \hol_j|_{\P(U_i \cap U_j)} 
  \]
  extending the transition function $g_{ij}$.
  In other words, for each path
  $\gamma \maps x \to y$ in $U_i \cap U_j$ a morphism in $\G$:
\vskip 0.5em
\begin{center}
\begin{picture}(100,100)
  \includegraphics{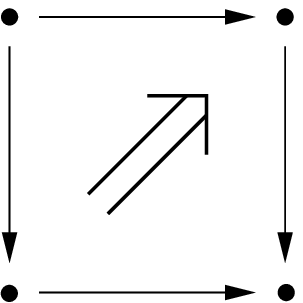}
  \put(-60,90){$g_{ij}\of{x}$}
  \put(-60,-8){$g_{ij}\of{y}$}
  \put(-120,42){$\hol_i(\gamma)$}
  \put(0,42){$\hol_j(\gamma)$}
  \put(-75,50){$g_{ij}\of{\gamma}$}
\end{picture}
\vskip 0.5em
\end{center}
\noindent
depending smoothly on $\gamma$, such that this diagram commutes:
\vskip 0.5em
\begin{eqnarray*}
\begin{picture}(240,140)
 \includegraphics{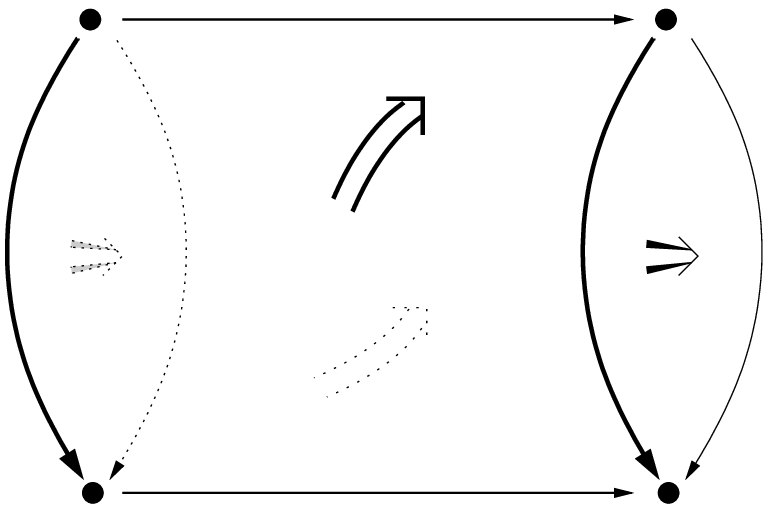}
 \put(-3,-17){	
 \begin{picture}(0,0)
 \put(122,0){
  \begin{picture}(0,0)
   \put(-385,90){$\hol_i\of{\gamma}$}
   \put(-290,100){${}_{\hol_i\of{\eta}}$}
   \put(-333,100){$\hol_i\of{\Sigma}$}
   \put(-242,50){${}_{g_{ij}\of{\eta}}$}
   \put(-268,125){${g_{ij}\of{\gamma}}$}
   \put(-250,163){$g_{ij}\of{x}$}
   \put(-250,10){$g_{ij}\of{y}$}
  \end{picture}
 }
 \put(288,0){
  \begin{picture}(0,0)
   \put(-382,90){$\hol_j\of{\gamma}$}
   \put(-290,100){${}_{\hol_j\of{\eta}}$}
   \put(-333,100){$\hol_j\of{\Sigma}$}
  \end{picture}
 }
 \end{picture}
}
\end{picture}
\end{eqnarray*}
for any bigon $\Sigma \maps \gamma \To \eta$ in $U_{ij}$,
\item
  for each $i,j,k \in I$ the transition function
  $h_{ijk}$ defines a modification:
\begin{center}
 \begin{picture}(100,80)
  \includegraphics{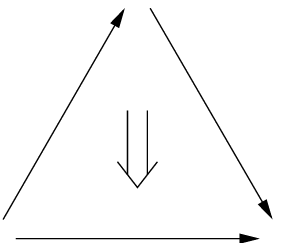}
  \put(-95,0){$\hol_i$}
  \put(-48,75){$\hol_j$}
  \put(2,0){$\hol_k$}
  \put(-78,40){$g_{ij}$}
  \put(-17,40){$g_{jk}$}
  \put(-50,-7){$g_{ik}$}
  \put(-64,24){$h_{ijk}$}
  \end{picture}
  \vskip 1em
\end{center}
In other words, this diagram commutes:
\begin{center}
\begin{picture}(180,270)
\includegraphics{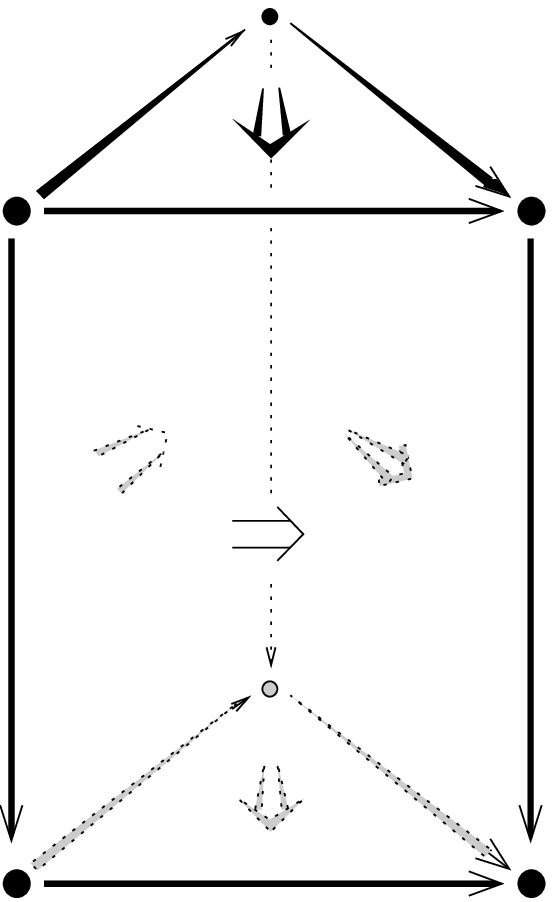}
  \put(-3,100){$\hol_k\of{\gamma}$}
  \put(-190,100){$\hol_i\of{\gamma}$}
  \put(-80,160){${}_{\hol_j\of{\gamma}}$}
  \put(-100,-7){$g_{ik}\of{y}$}
  \put(-100,187){$g_{ik}\of{x}$}
  \put(-138,43){${}_{g_{ij}\of{y}}$}
  \put(-138,237){${}_{g_{ij}\of{x}}$}
  \put(-47,43){${}_{g_{jk}\of{y}}$}
  \put(-47,237){${}_{g_{jk}\of{x}}$}
  \put(-95,16){${}_{h_{ijk}\of{y}}$}
  \put(-95,212){${}_{h_{ijk}\of{x}}$}
  \put(-100,90){$g_{ik}\of{\gamma}$}
  \put(-150,124){${}_{g_{ij}\of{\gamma}}$}
  \put(-50,117){${}_{g_{jk}\of{\gamma}}$}
\end{picture}
\end{center}
for any bigon $\Sigma \maps \gamma \To \eta$ in $U_{ijk}$.
\end{itemize}
\end{definition}

\noindent
In addition, it is worth noting that 
whenever we have a 2-connection, for each $i,j,k,l \in I$ 
this tetrahedron commutes:

\begin{eqnarray*}
\begin{picture}(140,140)
  \includegraphics{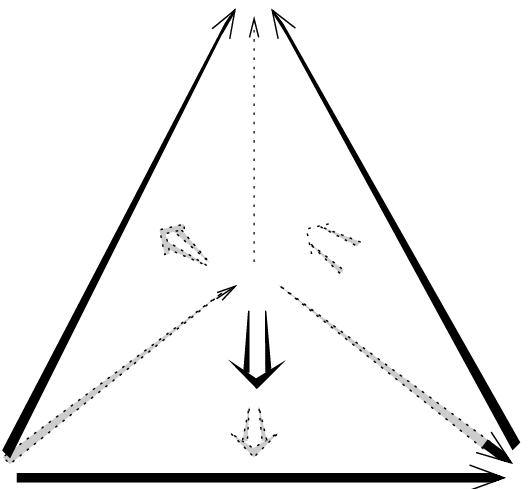}
  \put(-117,46){${}_{g_{ij}}$}
  \put(-50,46){${}_{g_{jk}}$}
  \put(-90,-7){$g_{ik}$}
  \put(-130,80){$g_{il}$}
  \put(-37,80){$g_{kl}$}
  \put(-78,101){${}_{g_{jl}}$}
  \put(-170,-3){$\hol_i$}
  \put(3,-3){$\hol_k$}
  \put(-85,62){${}_{\hol_j}$}
  \put(-85,145){${}_{\hol_l}$}

  \put(-100,14){${}^{h_{ijk}}$}
  \put(-70,30){$h_{ikl}$}
  \put(-68,85){${}_{h_{jkl}}$}
  \put(-110,80){${}_{h_{ijl}}$}
\end{picture}
\end{eqnarray*}

\vskip 1 em
In analogy to the situation for ordinary connections on bundles, one
would like to obtain 2-connections from Lie-algebra-valued
differential forms.  This is our next result.  In what follows,
$(G,H,t,\alpha)$ will be the smooth crossed module corresponding to
the smooth 2-group $\twogroup$.  We think of the transition function
$g_{ij}$ as taking values in $\Ob(\twogroup) = G$, and think of
$h_{ijk}$ as taking values in $H$.  Actually $h_{ijk}$ takes values in
$\Mor(\twogroup) \iso G \semidir H$, but its $G$ component is determined
by its source, so only its $H$ component is interesting.  In these
terms, the fact that
\[    h_{ijk}(x) \maps g_{ij}(x) g_{jk}(x) \stackto{\sim} g_{ik}(x)  \]
translates into the equation
\[  g_{ij}(x) \, g_{jk}(x) \, t(h_{ijk}) = g_{ik}(x)  ,\]
and the associative law of Theorem \ref{2-bundle cocycle conditions}
(i.e.\ the above tetrahedron) becomes a cocycle 
condition familiar from the theory of nonabelian gerbes:
\[
  h_{ijk} \, h_{ikl} 
 = 
  \alpha(g_{ij})(h_{jkl}) \, h_{ijl} 
 \,.
\]

\begin{theorem} \et
\label{central proposition concerning 2-connections}
   For any smooth 2-group $\twogroup$ and smooth space $B$,
   suppose that $E \to B$ is a principal $\twogroup$-2-bundle equipped with 
   local trivializations over open sets $\{U_i\}_{i \in I}$ covering $B$, 
   with the transition functions $g_{ij}$, $h_{ijk}$ and 
   $k_i$ given as in Theorem \ref{2-bundle cocycle conditions}.   
   Suppose for simplicity that $k_i = 1$.  
   Let $(G,H,t,\alpha)$ be the smooth crossed module corresponding to
   $\twogroup$, and let $(\g,\h,dt,d\alpha)$ be the corresponding 
   differential crossed module.
   Then there is a one-to-one correspondence between 
   2-connections on $E$ and Lie-algebra-valued 
   differential forms $(A_i,B_i,a_{ij})$ satisfying certain 
   equations, as follows:

\begin{itemize}

\item
  The local holonomy 2-functor $\hol_i$ is specified by 
  differential forms
  \begin{eqnarray*}
     A_i \in \Omega^1\of{U_i,\g}
     \\
     B_i \in \Omega^2\of{U_i,\h}
  \end{eqnarray*}
  satisfying
  \[
    F_{A_i} + dt\of{B_i} = 0
    \,,
  \]
  where $F_{A_i} = \extd A_i + A_i \wedge A_i$ 
  is the curvature 2-form of $A_i$.

 \item
   The transition pseudonatural isomorphism  
  $\hol_i \stackto{g_{ij}} \hol_j$
  is specified by the transition functions $g_{ij}$ together
  with differential forms
   \begin{eqnarray*}
     a_{ij} &\in& \Omega^1\of{U_{ij}, \h}
   \end{eqnarray*}
   satisfying the equations:
  \begin{eqnarray*}
    A_i &=& g_{ij}A_jg_{ij}^{-1} + g_{ij} \extd g_{ij}^{-1} - dt(a_{ij})
    \nonumber\\
    B_i &=& \alpha\of{g_{ij}}\of{B_j} + k_{ij}
  \end{eqnarray*}
  on $U_{ij}$,
  where
  \[
    k_{ij} = \extd a_{ij} + a_{ij}\wedge a_{ij} + d\alpha\of{A_i}\wedge a_{ij}
    \,.
  \]

 \item
  The modification 
  $g_{ij}\circ g_{jk} \stackto{h_{ijk}} g_{ik}$ 
  is specified by the transition functions $h_{ijk}$.  
  For this, the differential forms $a_{ij}$ are
  required to satisfy the equation:
  \[
  a_{ij} + \alpha(g_{ij}) a_{jk} = 
  h_{ijk} a_{ik} h_{ijk}^{-1} + (\extd h_{ijk})  h_{ijk}^{-1} + 
  d\alpha(A_i)(h_{ijk}) \, h_{ijk}^{-1} 
   \]
on $U_{ijk}$.

\end{itemize}
\end{theorem}

\Proof 
See Baez and Schreiber \cite{BaezSchreiber:2004}.
The `vanishing fake curvature' condition $F_{A_i} + dt\of{B_i} = 0$ 
is necessary for the holonomy 2-functor to preserve the source and
target of 2-morphisms.  It also guarantees that the holonomy over a
parametrized bigon is invariant under thin homotopies. 
\endofproof

The reader familiar with gerbes will recognize that Lie-algebra-valued
differential forms of the above sort give a \emph{connection on an
abelian gerbe} when $\twogroup$ is of the special form described in
Example \ref{crossed.1}.  Similarly, they give rise to a \emph{connection
with vanishing fake curvature on a nonabelian gerbe} when $\twogroup$
is of the form described in Example \ref{crossed.2}.

The vanishing fake curvature condition is a strong one.  As Breen has
emphasized, it implies that the $\h$-valued `curvature' 3-form $H =
\extd B + d\alpha(A) \wedge B$ actually takes values in the kernel of
$dt$, which is an abelian ideal of $\h$.  So, the existence of
well-behaved holonomies forces a 2-connection to be somewhat abelian
in nature.

\subsubsection*{Acknowledgements}

We are grateful to
Orlando Alvarez,
Paolo Aschieri,
Toby Bartels,
Larry Breen,
Dan Christensen,
James Dolan,
Jens Fjelstad,
Ezra Getzler,
Branislav Jur{\v c}o,
Mikhail Kapranov,
Anders Kock,
Amitabha Lahiri,
Thomas Larsson,
Hendryk Pfeiffer,
and 
Danny Stevenson
for comments and helpful discussion. 
The first author thanks the Streetfest organizers for inviting
him to speak on higher gauge theory along with Alissa Crans and
Danny Stevenson \cite{BCS}.
The second author was supported by SFB/TR 12.

\small{

}

\end{document}